# Some Observations on the Fig-8 Solution
# Off the 3-Body Problem

*F.Janssens, Feb 2004*
*Update Maart 2008*

## 1.Introduction

We discuss some properties of the periodic solution of the three-body problem where three particles of equal mass follow the same trajectory. This trajectory has the shape of a figure-8. The three particles have a constant separation in time of 1/3 of the period. The angular momentum of the system is zero. The paper by Broucke and Elipe[1] gives a thorough overview of history of this solution and places this solution in the general framework of the periodic solutions of the three-body problem. The two most interesting properties are the triangular and collinear configuration. Their occurrence in well known quartic curves with a figure eight shape are collected as an appendix.

## 2.Trajectory

In the inertial barycentric frame, the trajectory described by the 3 particles is shown in Fig.1. The fixed origin is the center of mass (COM) of three points with equal mass.

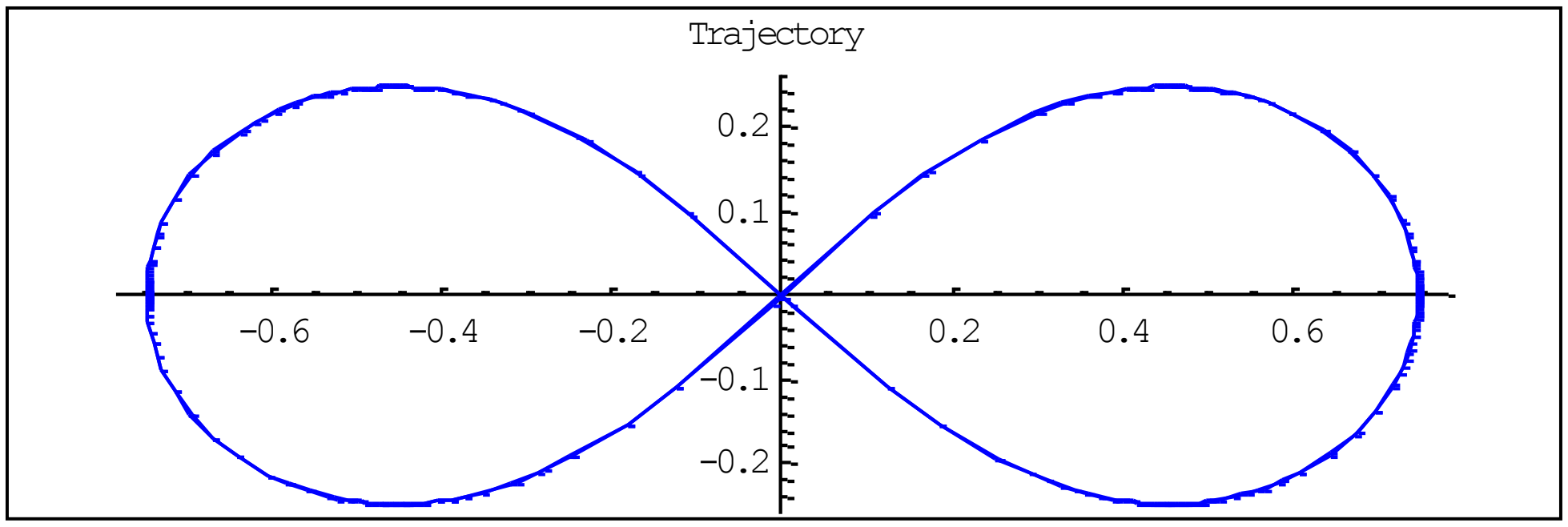

**Fig.1 – Trajectory figure-8**

### *1.Initial Configuration*

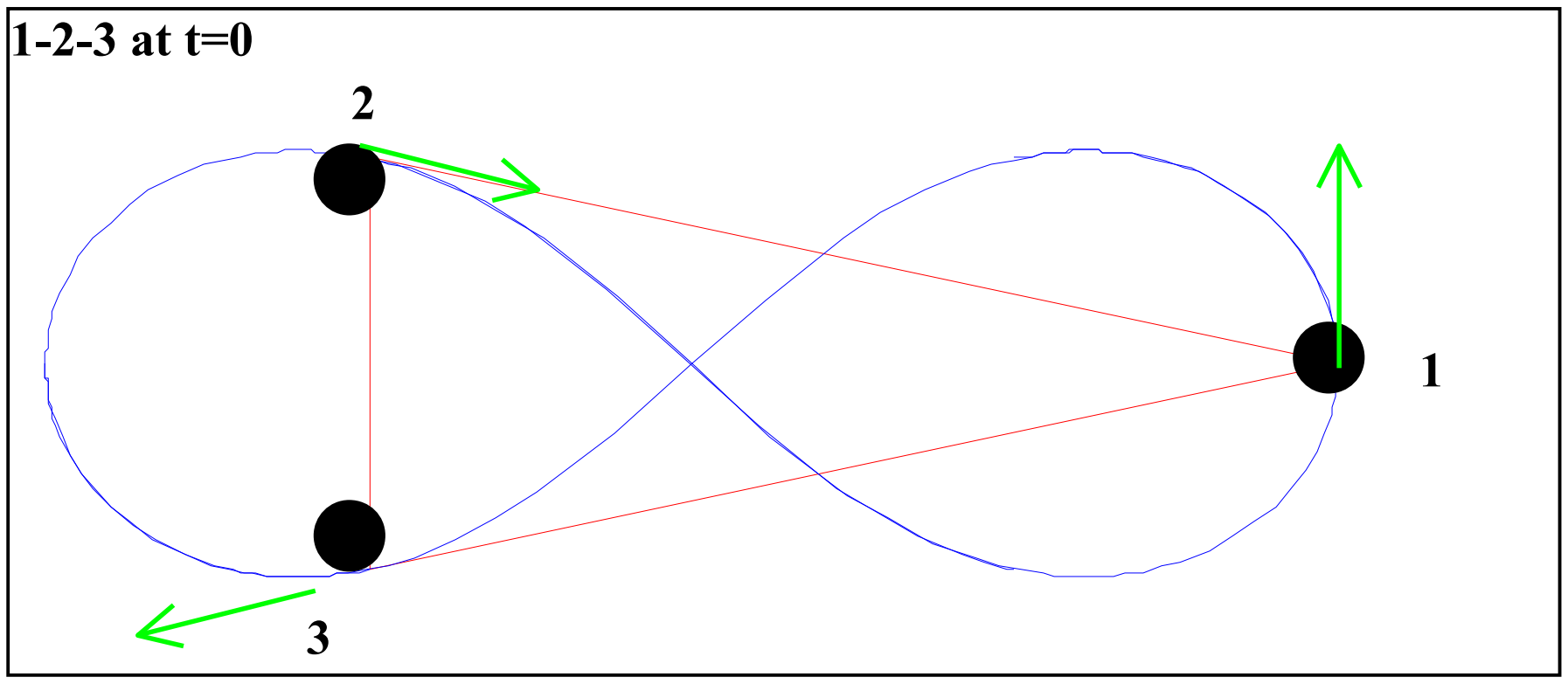

**Fig.2 – Initial Configuration – Isoceles Triangle**

As starting configuration for the 3 points, we take the point labeled 1 *at its maximal distance from the origin* and take that direction as x-axis. Such a point must necessarily exist for a periodic, hence bounded, solution. Its coordinates are ($x_0$ ,0) and the velocity is perpendicular to the x-axis. The direction of the velocity defines the y-axis and the sense in which point 1 and the two others run around the trajectory. From the well known first integrals, we have some information on the corresponding location for the two other points. As the COM is the fixed origin, we have at all t :

$$\bar{r}_1 + \bar{r}_2 + \bar{r}_3 = 0 \qquad (1)$$

because the three masses are equal. The conservation of linear momentum is given by the derivative of eq.(1) .The second derivative relates the accelerations $(\bar{a}_i)$ of the 3 points.

$$\begin{aligned} \bar{v}_1 + \bar{v}_2 + \bar{v}_3 &= 0 \\ \bar{a}_1 + \bar{a}_2 + \bar{a}_3 &= 0 \end{aligned} \qquad (2\&3)$$

From Eq.1 follows immediately that the COM of points 2 and 3 is at [-$x_0$/2,0] and the initial conditions of point 2 and 3 can be written as :

$$2 \Rightarrow [-\frac{x_0}{2} + \delta, y_0], \quad 3 \Rightarrow [-\frac{x_0}{2} - \delta, -y_0] \qquad (4)$$

For the initial velocities, we have a similar relation (Eq.2 ):

$$2 \Rightarrow [v_{x0}, -\frac{v_0}{2} + \delta v], \quad 3 \Rightarrow [-v_{x0}, -\frac{v_0}{2} - \delta v\ ] \qquad (5)$$

Using the fact that the total angular momentum is zero, gives :

$$v_{x0} y_0 = 3\,(v_0/2)\,(x_0/2) + \delta\ \delta v \qquad (6)$$

The simplest way to satisfy Eq.(6) giving maximal symmetry is $\delta = \delta v = 0$.
Points 2 and 3 have then the coordinates [-$x_0$/2, **±**$y_0$] and velocities [**±**$v_{x0}$, -$v_0$/2]. They are symmetric w.r.t. the x-axis. The three points make an isosceles triangle.

Rewriting Eq.(6) as:

$$\frac{y_0}{3x_0/2} = \frac{v_0/2}{v_{x0}} \qquad (7)$$

shows that the velocities have the same direction as the two sides of the isosceles triangle. The velocity of particle 2 points to particle 1 while the velocity of point 3 points in the opposite direction. Hence, the trajectory has the following property:

*The tangents from the points at maximal distance from the origin touch the curve at points on the other side of the origin at half the maximal distance on the x-axis.*

The non-dimensional numerical values given in Broucke[1] satisfy this property and are illustrated by the red triangle in fig.. We used them for the simulations used in this paper.

**R1** [ .746156, 0] ; **R2** [-0.373078, .238313]; **R3** [-0.373078, - .238313] (8)
**V1**[ 0, .324677] ; **V2** [.764226 , -.162339]; **V3** [-[.764226 , -.162339]

The sides of the triangle make an acute angle of ±12.05 deg. with the x-axis.

## *2. Normalization -Family of trajectories.*

The non-dimensional numerical values in Eq.(8) are normalized a reference length R , and total mass m. Each particle has mass m/3

- The lengths are fractions of R , the time is replaced by an angle $\tau = \sqrt{Gm/R^3}\ t$.
- The period in $\tau$ is $2\pi$ and $P = 2\pi\sqrt{R^3/Gm}$ in t.
- The velocities are in fractions of $\sqrt{Gm/R}$
- The accelerations are in units of $Gm/R^2$
- The angular momentum is in $\sqrt{GmR}$, the energy in Gm/R.

For example when the mass of the Earth is split over 3 particles or spheres of equal mass, and R equals the equatorial radius $R_E$ of the Earth, the period would be the period of a grazing satellite or 84.486 minutes. The maximal separation from the COM of the three masses would be $R_E$ x .74615 = 4758.9 km. The maximum and minimum separation between the 3 spheres is 7298.5 and 3039.94 km. With the density of earth, each sphere has a radius of 4671.86 km. The density of the three spheres must at least be 4 times the density of the Earth to avoid collisions.

For m = 1kg (one particle .333 kg) and R=1m, the period is ~28.153 days. Each particle is at most 74,6 cm from their COM. The separation varies between 1.12 and .47 m.

The normalization contains only one free parameter R when the total mass is fixed. The initial conditions have to be as above in Eq.8 or any other point on the fig-8 curve.

*This defines a one-dimensional line in the manifold of planar solutions in the phase space.* Initial conditions not on this line give, in general, not periodic solutions.

## *3. Description trajectory - Collinear Configuration.*

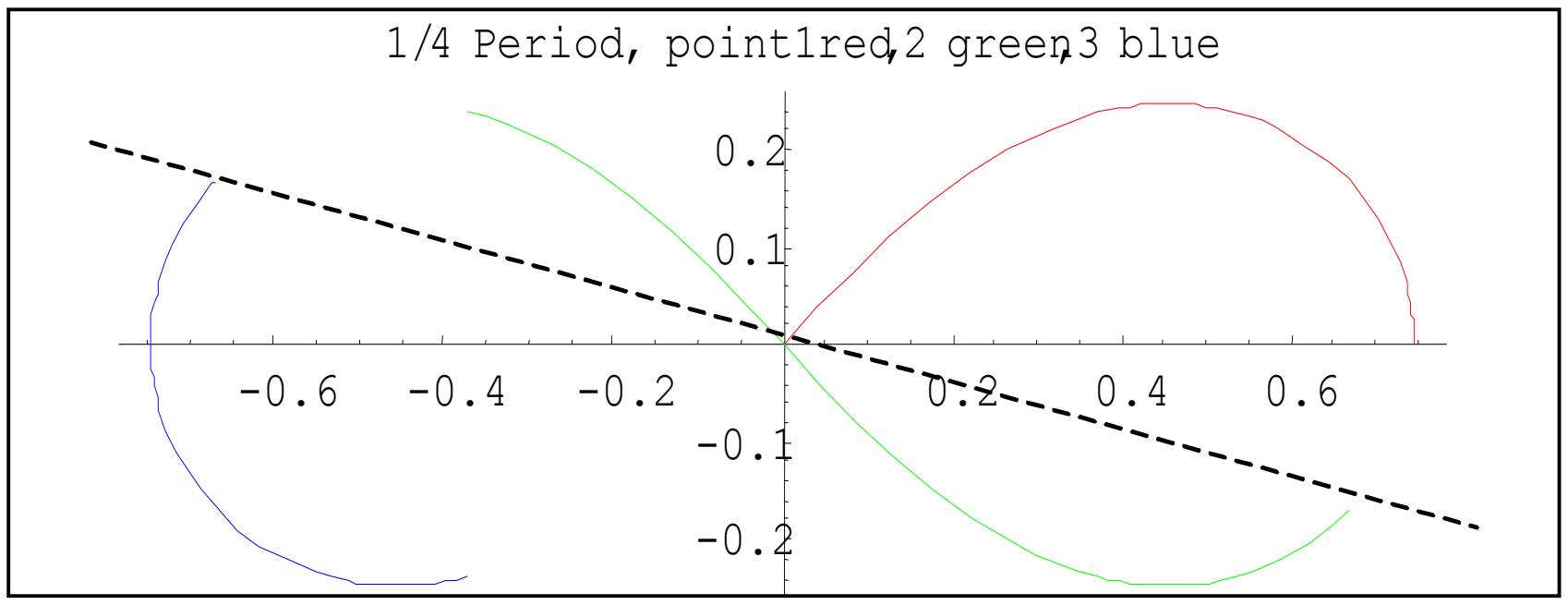

**Fig.3 – Trajectory over P/4.**

Starting from the isoceles configuration of Fig.2, we obtain fig.3 after a quarter of a period. Particle 1 is at the origin and did exactly . a quarter of the trajectory. Particle 2 has gone through the origin and is far in the right part of the eight while particle 3 is still in the left side of the eight.

When particle 1 is at the origin, then from the immobility of the COM, particles 2 and 3 are necessarily on a line passing through the origin and at the same distance from O (dotted line in Fig.3) . *Such collinear configuration* occurred also earlier when particle 2 was at the origin.
The integration showed that particle 2 was at the origin at $\tau$ = 30 deg. The position of the 2 other particles was :

**1**[.669531, .167788] , **3**[[-0.669531, -0.167788]

The two collinear configurations make an angle of ±14.0688 deg with the x-axis.

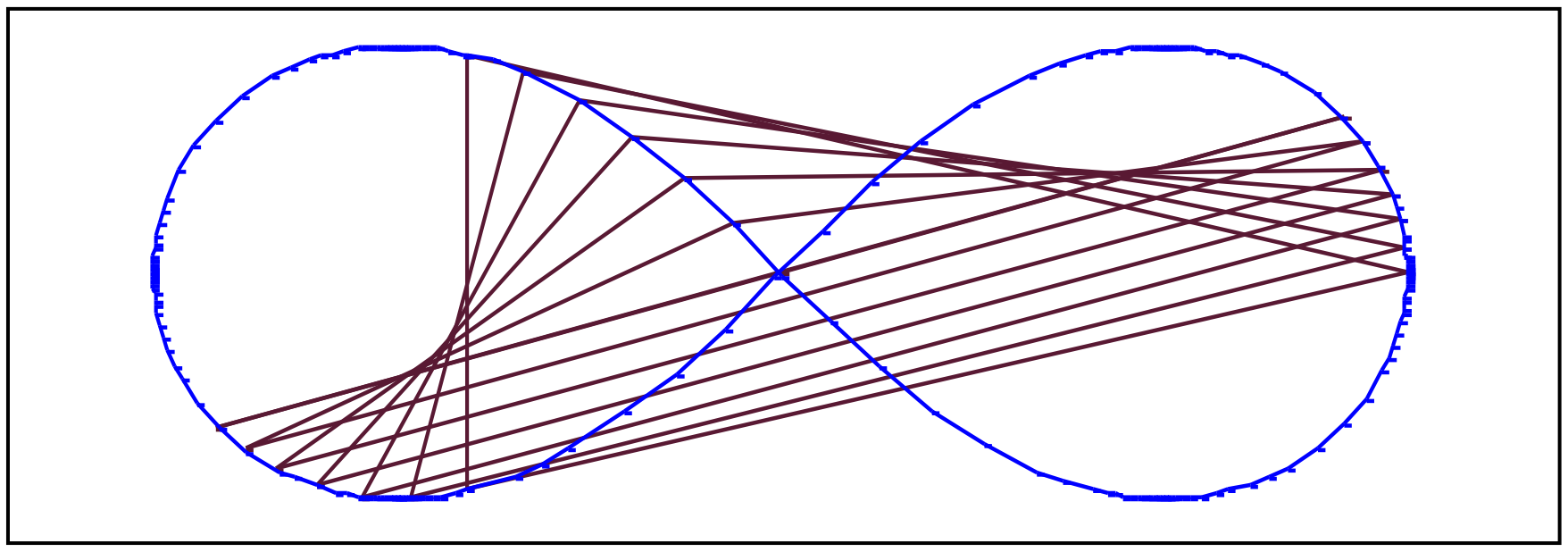

**Fig.4 – Evolution of triangle from $\tau$= 0, 30 deg**

Fig.4 shows the evolution from the iscoceles triangle at $\tau$ = 0 to the collinear configuration at t =30 deg. When $\tau$ increases further to 60 deg., we have again an isoceles configuration with particle 3 at the top and particle 1 and 2 at the base ( Fig.5) As the trajectory is symmetric w.r.t. the x- and y-axis, it is also symmetric w.r.t. the origin. This implies that the tangents (velocities) are parallel at points with coordinates [x,y] and [-x,-y] . This is the case for the two particles not at O, in the collinear configuration. As the angular momentum of the particle at O is zero, the angular momentum of the 2 other particles must cancel out which implies now that their velocities are equal and not only the components normal to the radius vector. ).

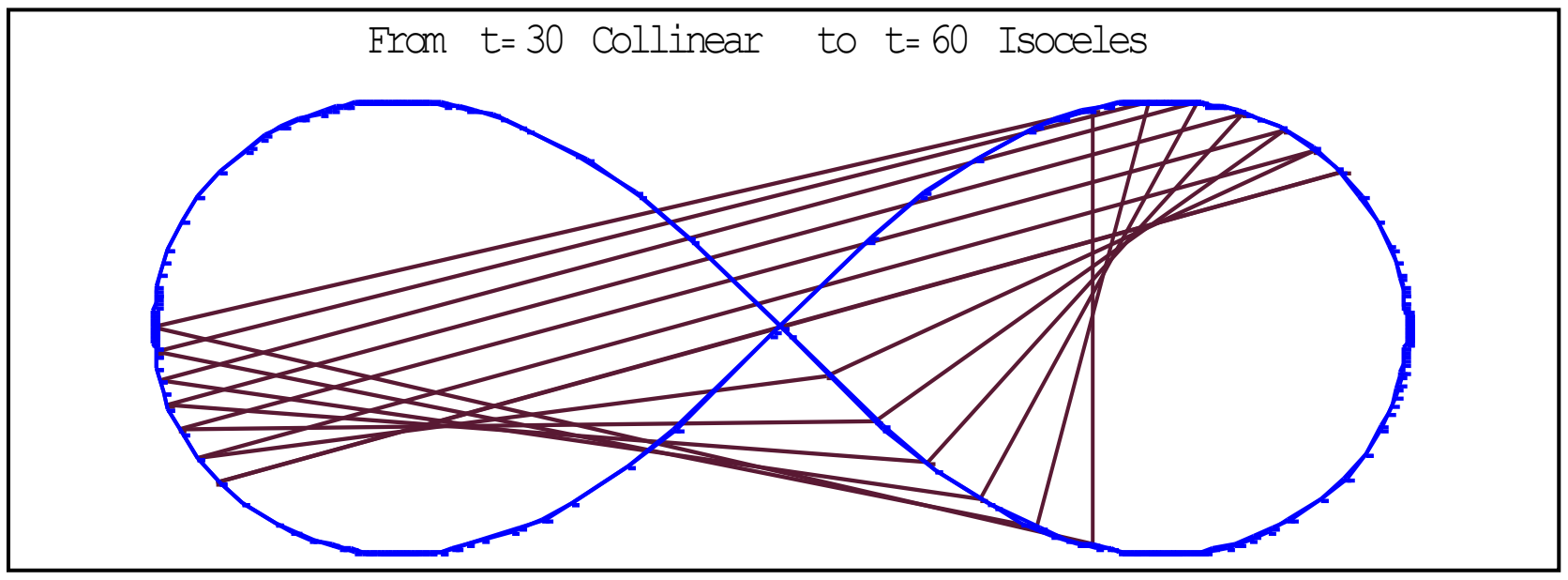

**Fig.5 – Evolution of triangle from from $\tau$=30, 60 deg**

*In a collinear configuration, the three tangents to the trajectory are parallel.* This direction makes an angle of ± 42.84 deg with the x-axis. Further more, as the sum of the 3 velocities is zero, *the velocity of the particle at the origin equals twice the velocity of any of the others with the direction reversed (Eq.2)*. The acceleration on the particle at the O is zero. The velocity of the particle at the origin is constant and is maximal (see Fig.x).

Now, a closed trajectory has necessarily a minimum and maximum value for the distance r from the origin. These two values are normally found by solving $\dot{r}$ =0. The particle at O has r = 0 which is clearly the minimal value but $\dot{r}$ is different from zero. The only points of the trajectory where $\dot{r}$ =0 is for the maximal value of r.
Fig.6 gives the 12 successive passages through an isoceles and collinear configuration which occur every 30 deg , during a period:

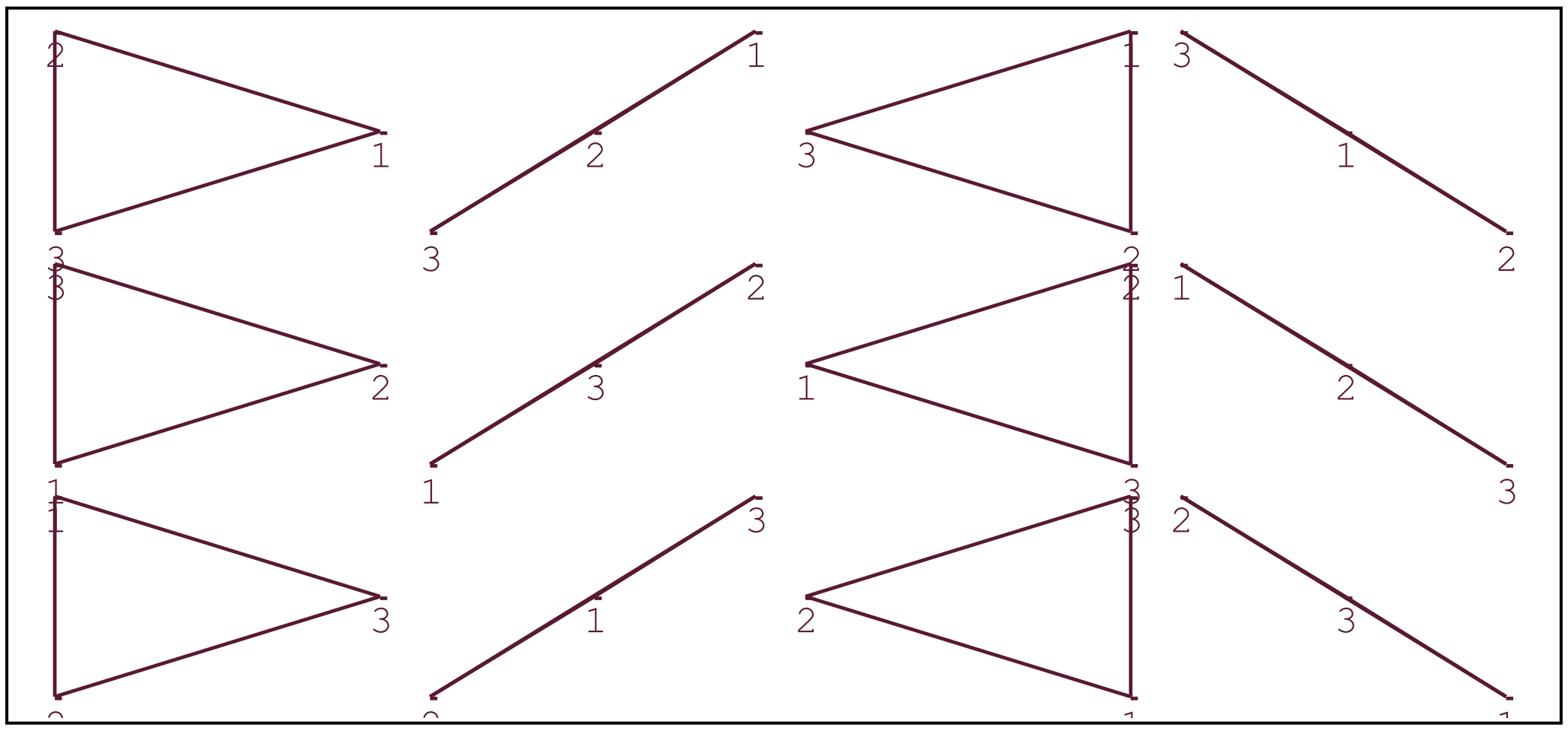

**Fig.6- Sequence every 30 deg. from left to right and top to bottom.**

The collinear configurations occur every 60 deg. At those instants there is a switchover for the loop of the figure-eight trajectory that contains two particles. The lines joining the two particles in the same loop define an envelope as shown in figures 7-8.

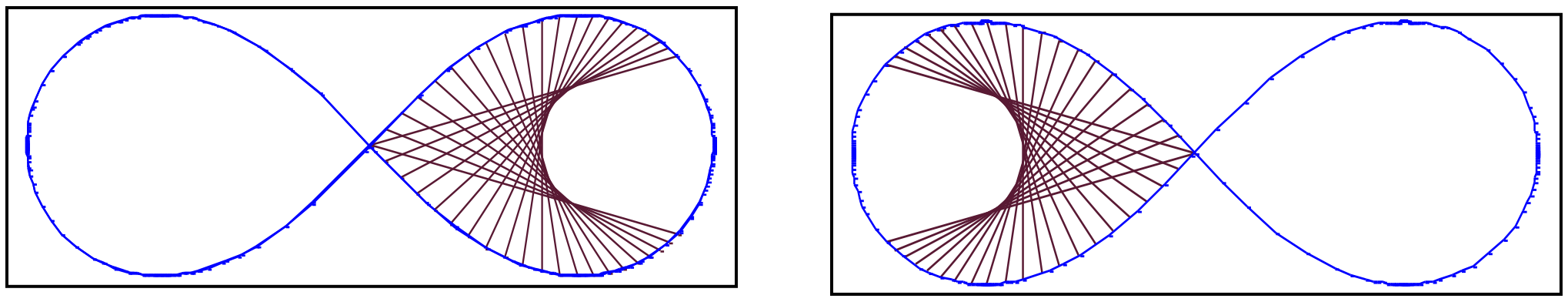

**Fig.7-8 Envelope defined by particles in the same loop**

Between two collinear configurations, the three particles together describe one loop of the curve. For instance, when $\tau$ goes from 30 to 90, particle 3 describes the part of the figure-8 with x-coordinate smaller as -.669531. Particle 2 starts from the origin and describes the lower part of the curve to x=.669531. Particle 3 starts with this x-coordinate, describes the upper part of the curve and ends at the origin.

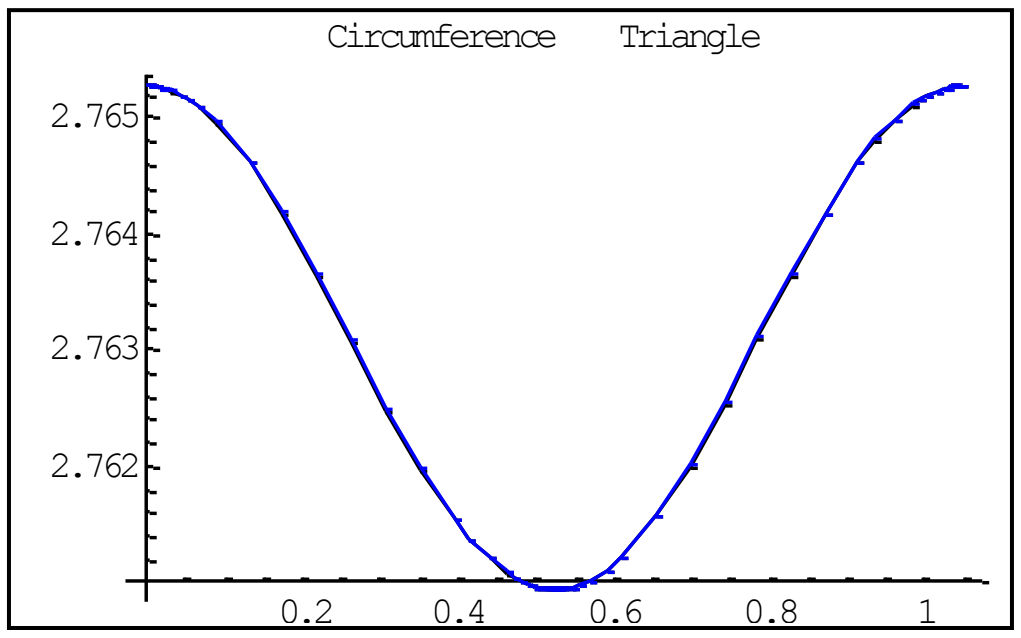

**Fig9. - Circumference of the triangle**

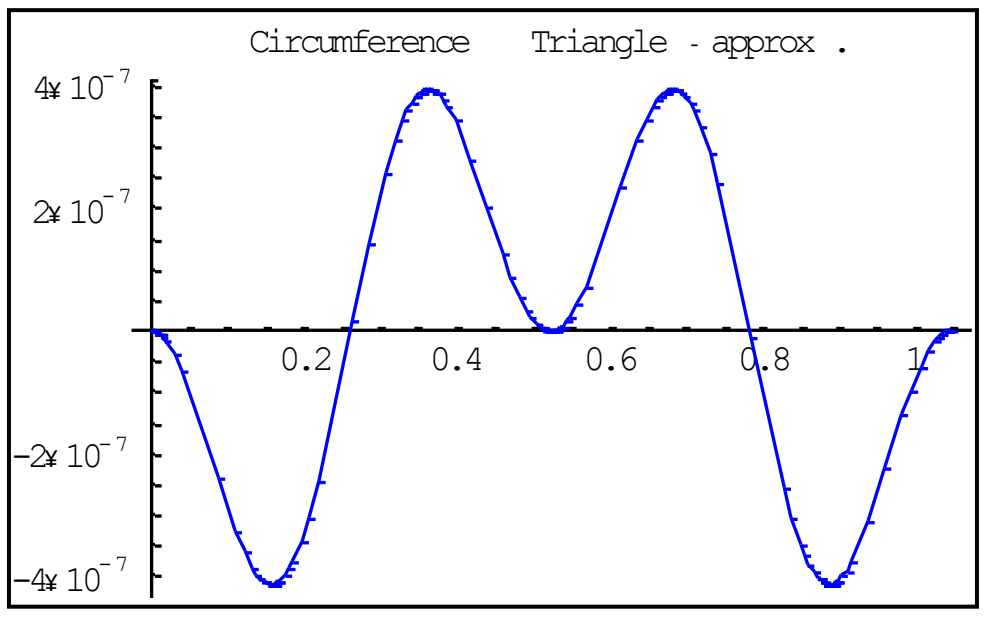

**Fig.10 - Approximation**

The triangle evolves from an isoceles to a degenerate collinear configuration . As the particles never collide, none of its sides becomes ever zero.
The evolution of the length of one side of the triangle is shown in Fig.12, for each side.

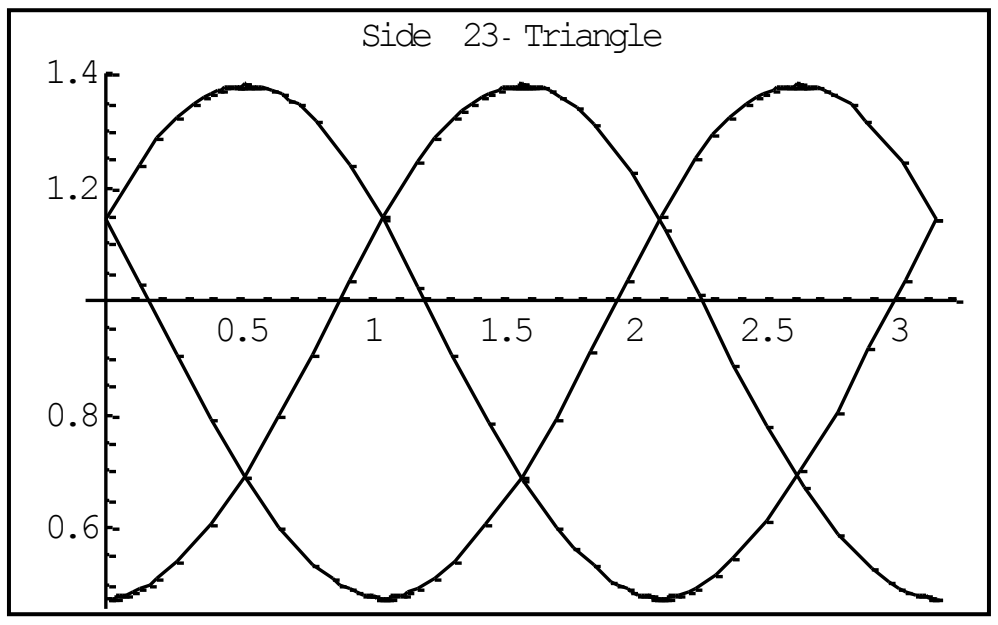

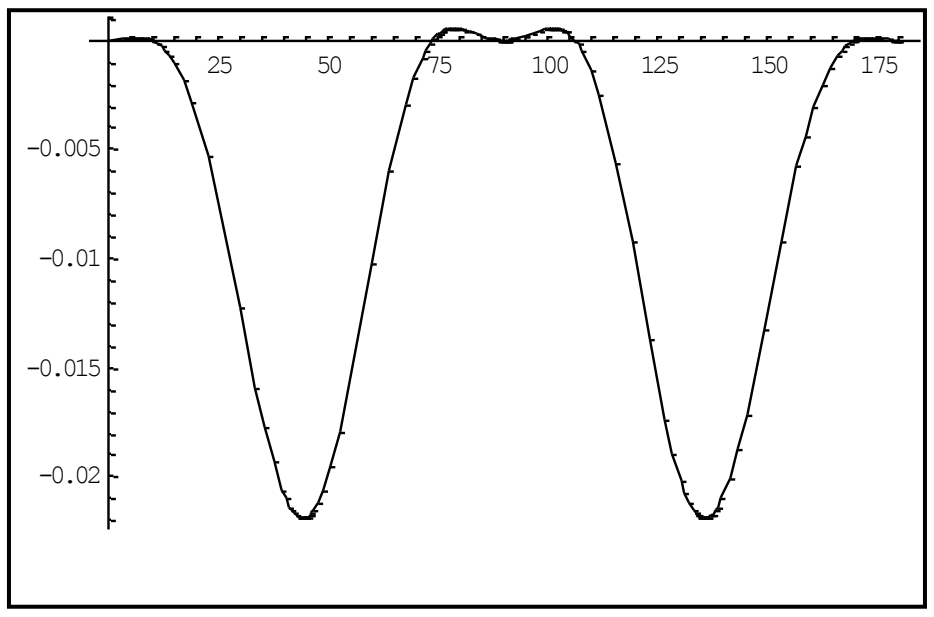

**Fig.12 – Separation $r_{ij}$** **Fig.13 – Residu between Circumference and Approximation**

The maximal length occurs in the collinear configuration between the 2 particles at the endpoints of the segment: $\ell_M$ =1.38047.The minimum length is the base of the iscoceles triangle $\ell_m$ =. 476626.

A first approximation , e.g. for $r_{23}$ (the curve starting from its minimum) , is given by

$$r_{23}(\tau) = .928549 - .451923\cos 2\tau$$

A good approximation of the circumference s is :

$$s(\tau) = c_0 + c_1 \cos 6\tau + c_2 \cos 12t$$

$c_0$= 2.76309294639 ; $c_1$=.00216672; $c_2$=.0000158567
The residu is shown in Fig.13.

*4. Velocity*

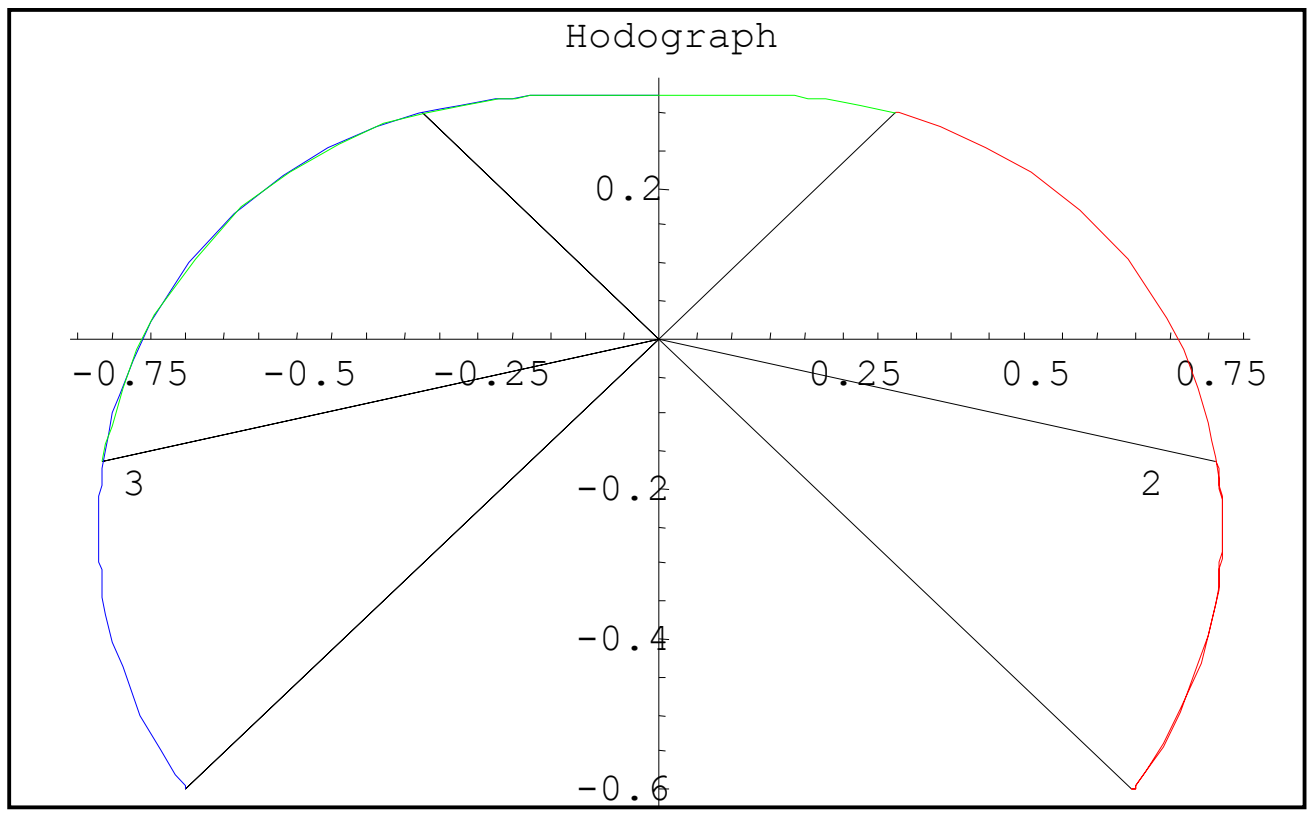

**Fig.14 – Hodograph or velocity diagram**

The hodograph is shown in Fig.14. At τ = 0, the velocity of particle 1 on the y-axis (vertical) and starts turning to the left as the acceleration is perpendicular to it. Particles 2 and 3 are the base of the triangular configuration in the left loop (Fig.7). The endpoints of $V_2$ and $V_3$ at τ = 0 shown in Fig.14. $V_2$ starts moving down, while particle 3 moves up. When particle 2 is at the lowest point (τ = 30), where the velocity is maximal, $V_1$ and $V_3$ meet in the opposite direction. The velocity of particle 2 makes then an angle of –42.84 deg. with the x-axis. The velocity of particle 1 is horizontal for $\tau$ = 53.66 deg. In half a period, $V_1$ has gone to its lowest point (222.84 deg.) and returned to its initial orientation of 90 deg. The total rotation of the velocity vector of a particle over half a period is zero. *The hodograph is not a closed curve.* Each particle completes the curve in a period with a total rotation of zero.

## 5.Acceleration

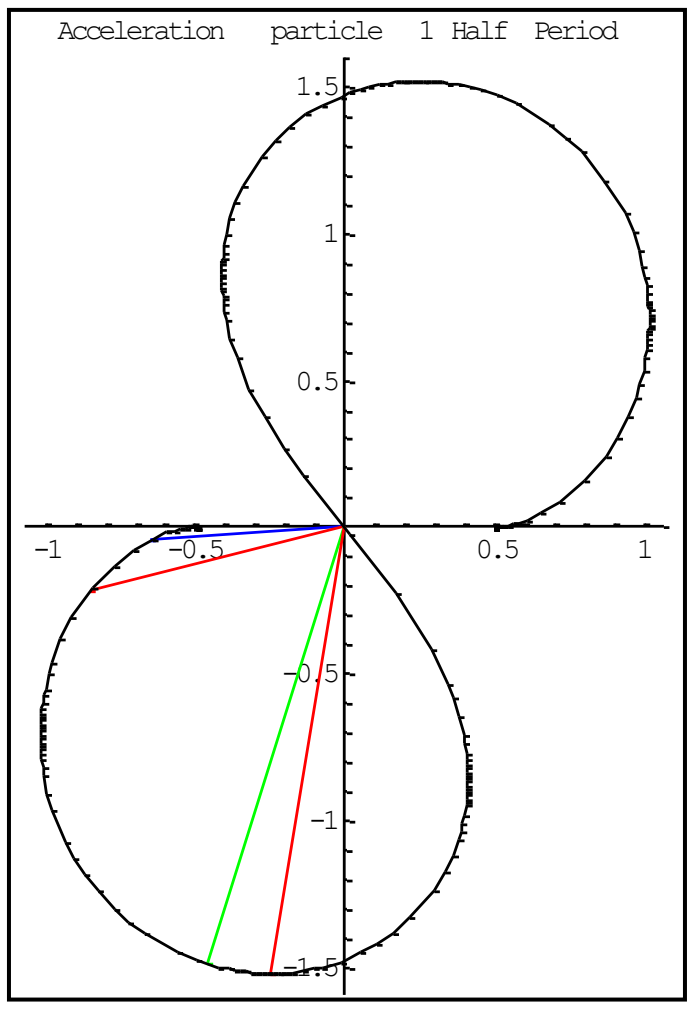

**Fig.15 - Acc. particle1 [0 - π ]**

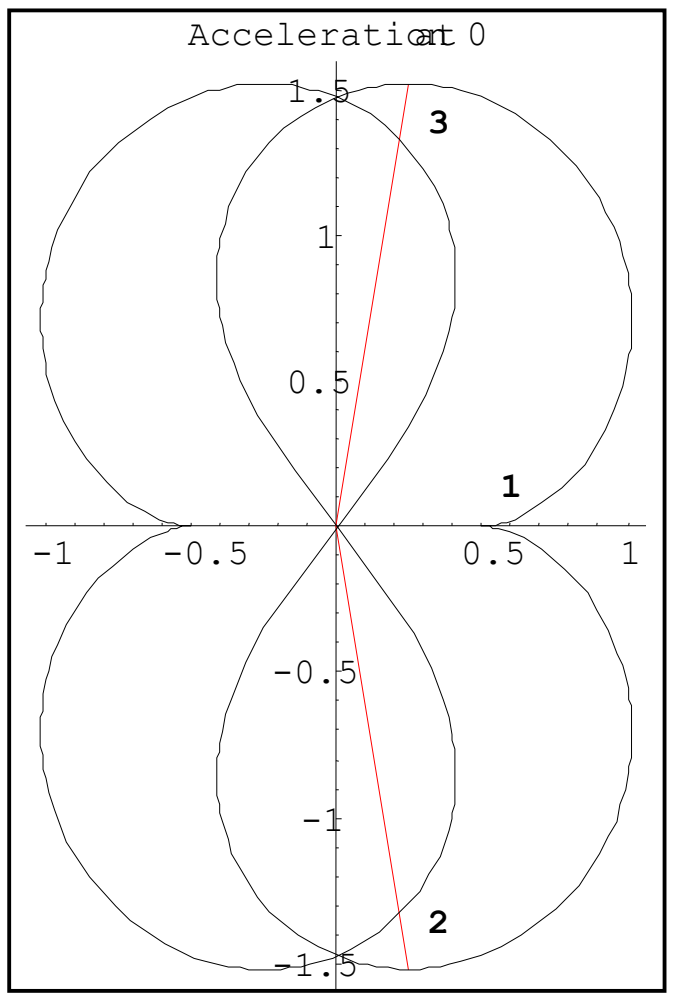

**Fig.16 - Acceleration [0=2 π ]**

Fig.15 shows the acceleration diagram of particle 1 over the first half period.

<u>$\tau$ = 0</u>, the acceleration is horizontal and pointing to the left as the configuration of the three particles is an isoceles triangle (Fig.1). Then, as particle 2 approaches particle 1 from above and particle 3 recedes while below particle 1, the acceleration of particle 1 should point slightly upwards. This is not visible on Fig.15, but zooming in for small values of t, (Fig.17), we see that this is indeed the case for a short time.

<u>$\tau$ = 7.19</u>, the acceleration is again horizontal and continues to rotate downwards. During the first quarter period, particle 2 gives the more important contribution to the acceleration of particle 1 as it is the closest.

<u>$\tau$ = 19.286</u>, particle 1 and 2 has the same y-coordinate (1 going up and 2 going down).The corresponding acceleration is the blue line in Fig.15.

When <u>$\tau$ =30</u>, the configuration is collinear. The acceleration points to origin where particle 2 is at that instant. The angle with the x-axis is 14.0688 deg (red line in Fig.16). Particle one continues to move upwards till <u>$\tau$ = 53.663</u>.

<u>$\tau$=57.19</u>, just before the isosceles configuration with 3 as top and 1-2 as base, the acceleration of 1 is maximum (green line in Fig.15).

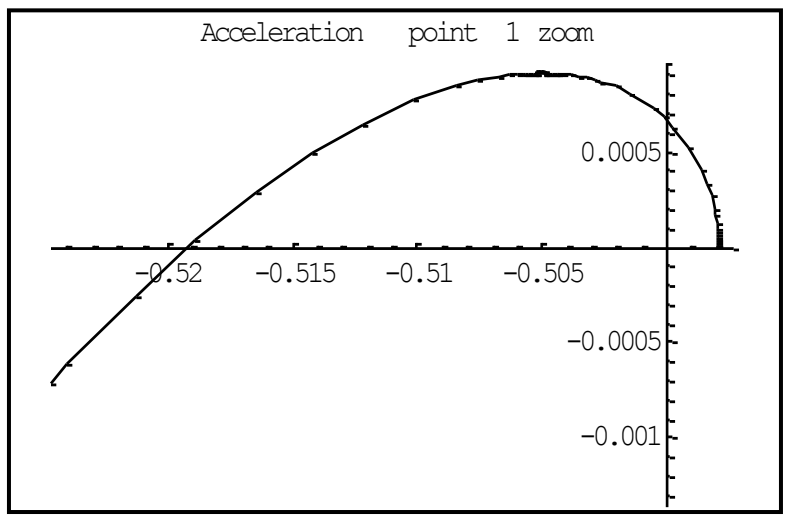

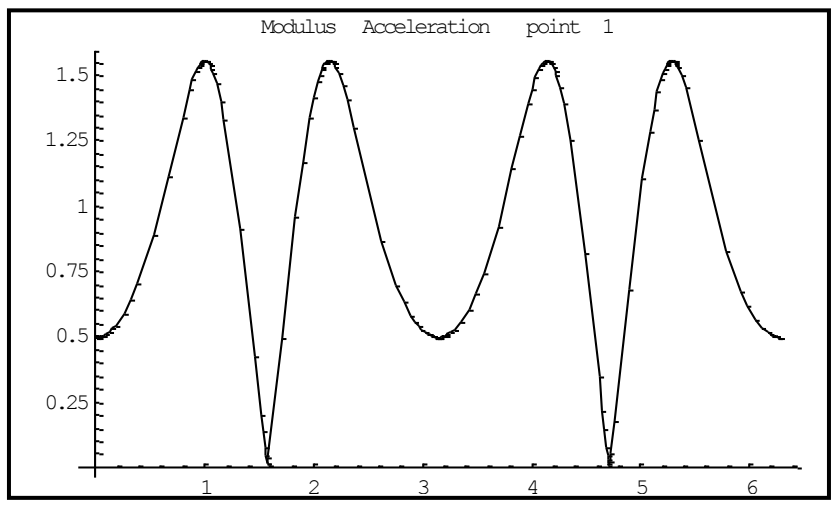

**Fig.17- Zoom accel.1 $\tau$= 0-8 deg** **Fig.18 - Modulus acceleration**

$\underline{\tau = 60}$ (isosceles configuration ) the acceleration of particle 1 is the second red line in Fig.15. From that instant on, particle 2 is to the right of particle 1 but the acceleration is still to the left due to particle 3.

$\underline{\tau = 63.431}$ , the acceleration is vertical downwards and will be to the right from then onwards. At a quarter period, $\underline{\tau = 90}$, particle 1 is the origin, its acceleration is zero and the configuration again collinear. When particle 1 approaches the origin, the direction of its acceleration goes to –51. 5693 deg.

For $\tau > 90$, the acceleration rotates back and is the symmetric one w.r.t. the origin of the point 180-$\tau$. Point 1 is back on the x-axis, after half a period or $\tau$=180.

Fig.16 shows the acceleration over a full period. This curve is, of course, also valid for particle 2 and 3 that start at a different point. The total rotation of the acceleration of each particle over a period is zero. The two inner lobes, symmetric to the y-axis are very clear on Fig.16. There are also two such small lobes , symmetric to the x-axis where the outer lobes seem to touch the x-axis. Fig.18 shows the modulus of the accelerations.

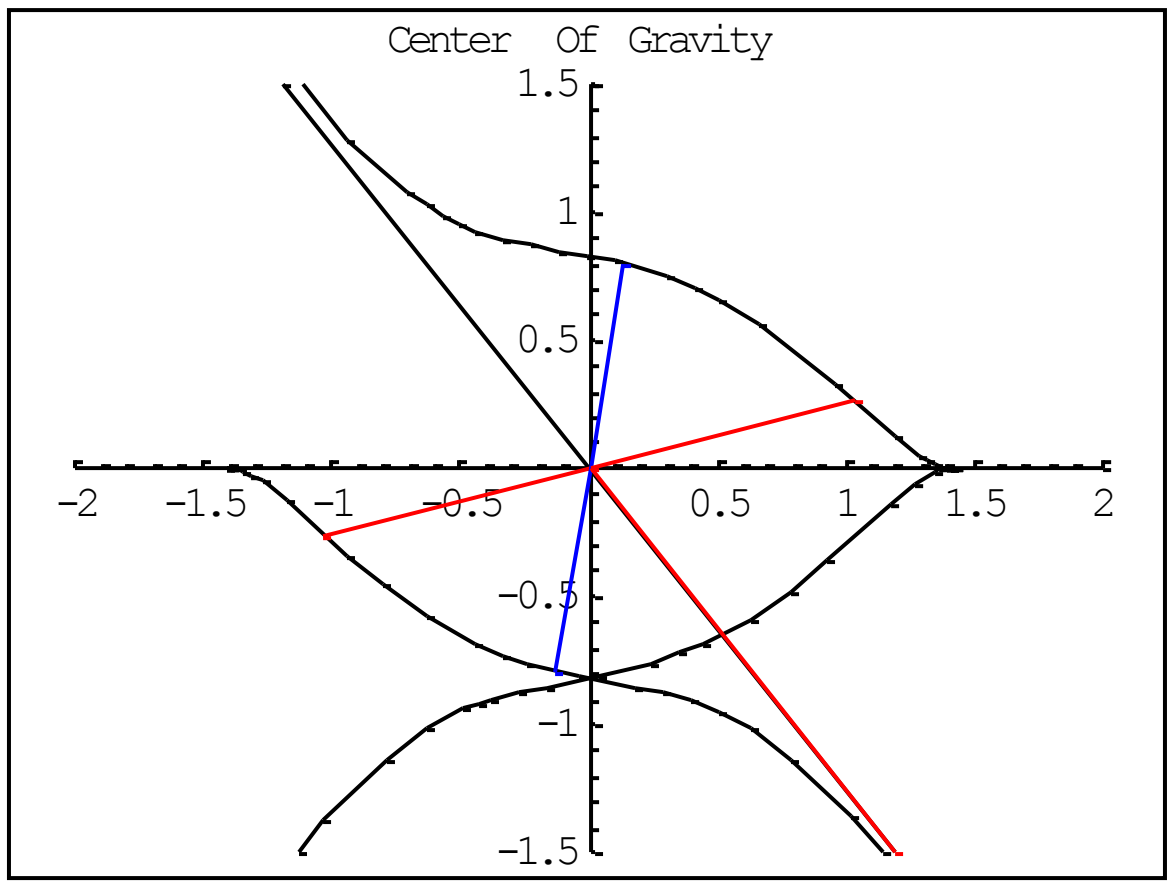

**Fig.19 – Center of Gravity**

These results for the acceleration can also be formulated in terms of the center of gravity (COG). This is the position of a particle with twice the mass of a single particle causing the same acceleration as two particles on the remaining one. This curve is shown in Fig.20 and is the inverse of the previous acceleration diagram. At $\tau$ =0, the COG is on the minus x-axis. The red and blue lines, in quadrant 3, indicate the successive collinear and triangular configurations. When particle 1 is at the origin, its acceleration is zero and the COG is at infinity. The direction of the asymptote is -51.5693 deg. The COG jumps from quadrant 4 to quadrant 2, moves to the next triangular configuration in quadrant 1 to continue to the x-axis ( triangular configuration after half a period ).

The second half period, the COG describes the symmetric curve w.r.t. the y-axis.

## *6.Relative Motion*

The two body problem is solved by taking one particle as the origin of an axis system while keeping the direction of the axes inertially fixed. The trajectory of the second particle is then found in this axis system. In inertial space (origin at the COM) both particles have a similar trajectory. Fig.21 shows the trajectory of particle 2 and 3 in such axis system centered at particle 1. These two particles follow a different trajectory that is not closed. One trajectory is the symmetric one of the other w.r.t. the origin.

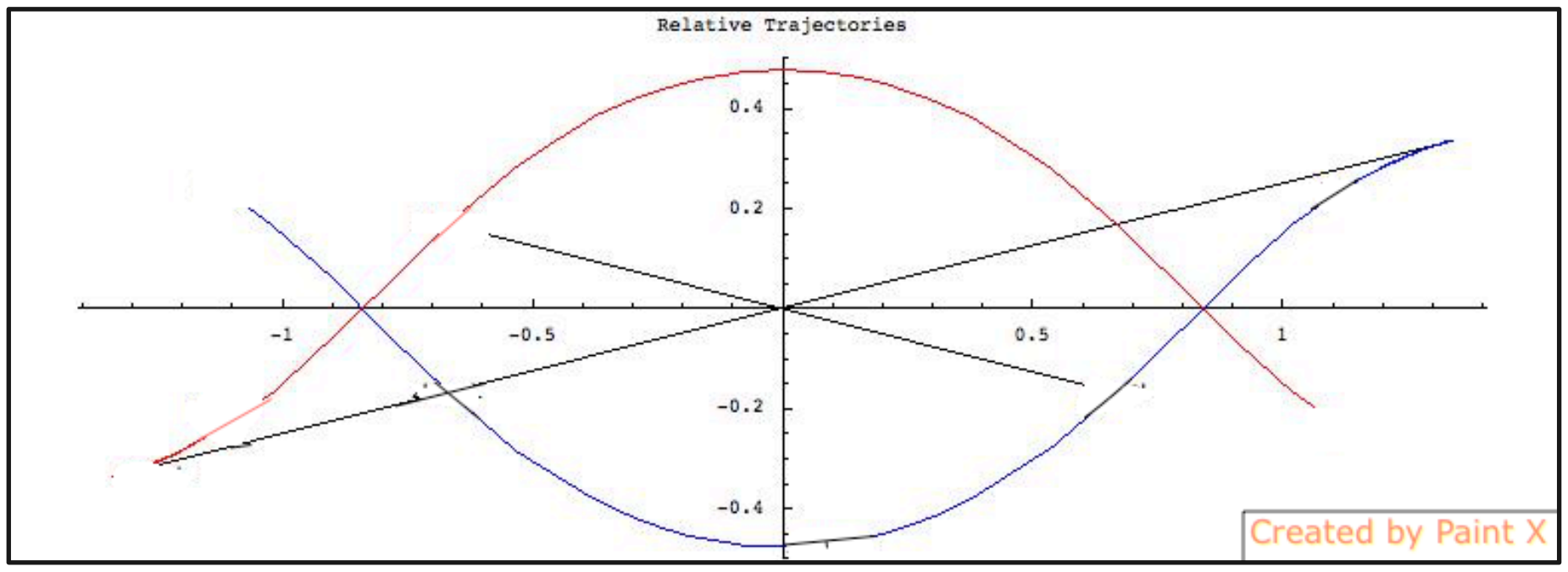

**Fig.20 – Relative trajectory of particle 2,3 w.r.t. particle 1**.

The particles run forwards and backwards on these curves. After one period they have returned to their starting point. However, the relative motion of the midpoint of the line joining particle 2 and 3, or their COM, describes a figure-8 curve similar to curve followed by the 3 particles in inertial space (size 3/2). In the inertial space, the COM of 2 and 3 follows the same fig-8 of half the size as the curve followed by the three particles. These statements follow immediately because the COM of the system is at rest.

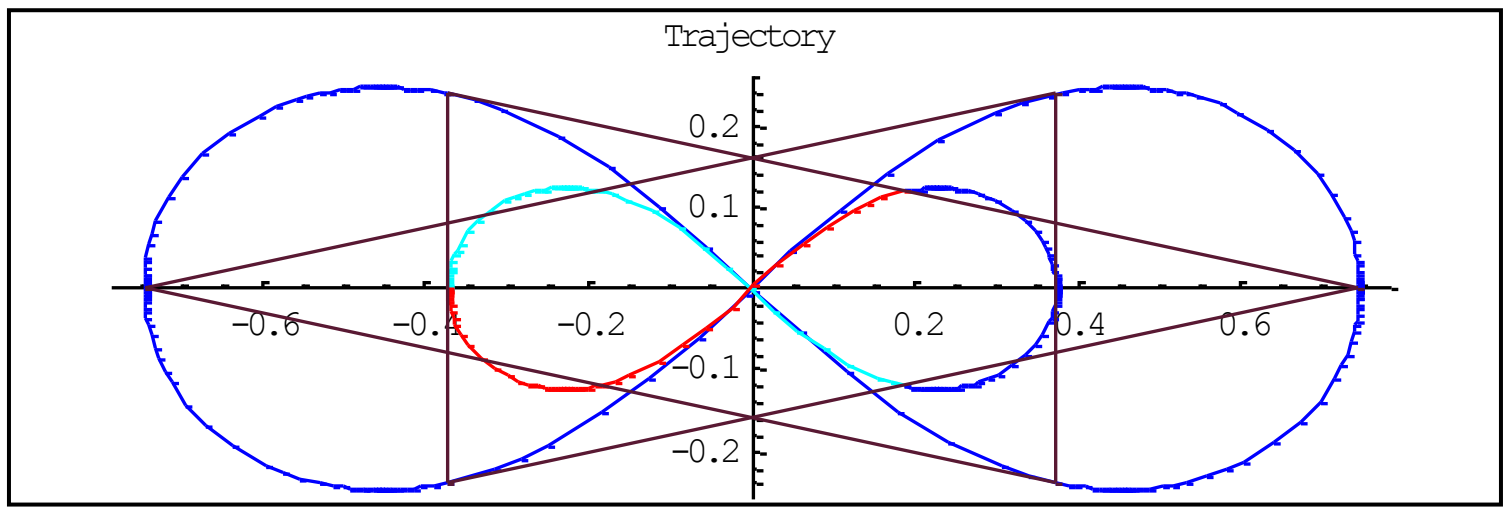

**Fig.21 – Trajectory of the COM of 2 and 3**.

## *7.Energy and angular momentum*

The total angular momentum is zero. At each instant, the angular momentum of one of the particles compensates the angular momentum of the two other ones. See Fig.22. As a consequence the average angular momentum of a single particle over a period is zero.

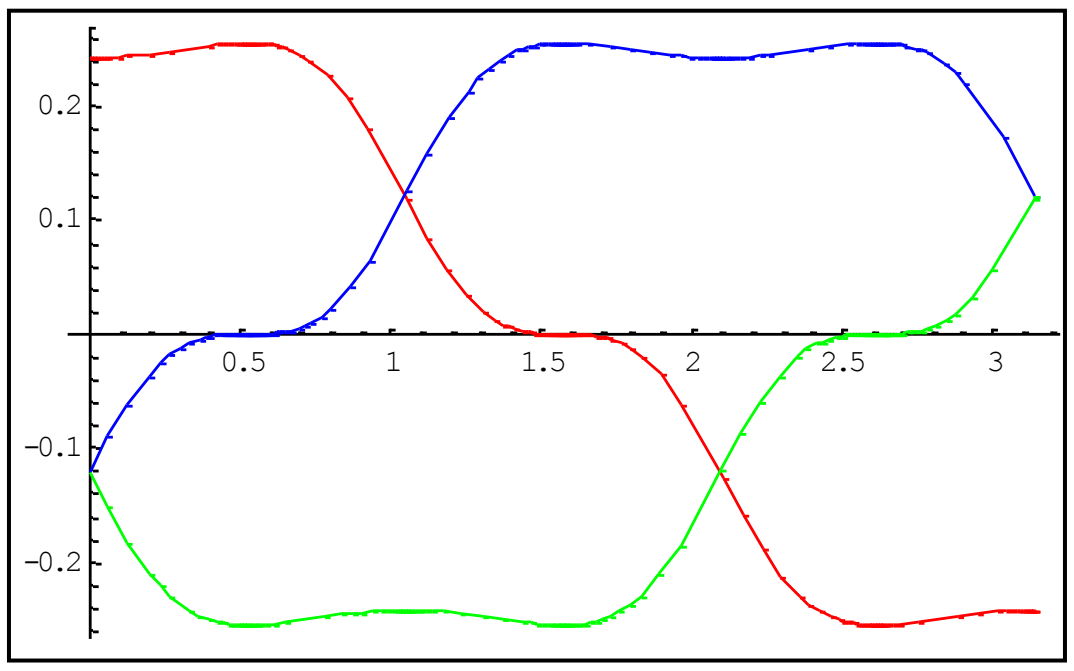

**Fig.22 – Angular Momentum of the three particles**

The energy w.r.t. inertial space : $E_{k3} = \frac{1}{2}\sum_{i=1,3} m_i V_i^2 \quad E_{p3} = \sum_{i \neq j=1,3} \frac{Gm_i m_j}{|r_i - r_j|}$ (9)

Fig.23 and 24 show the variability of the kinetic energy of one particle and inside a period and the potential energy of the system. The total energy of the trajectory is :

$E_t = - 0.6215955544$

In the two-body problem ,the most used expression is the relative energy per unit mass :

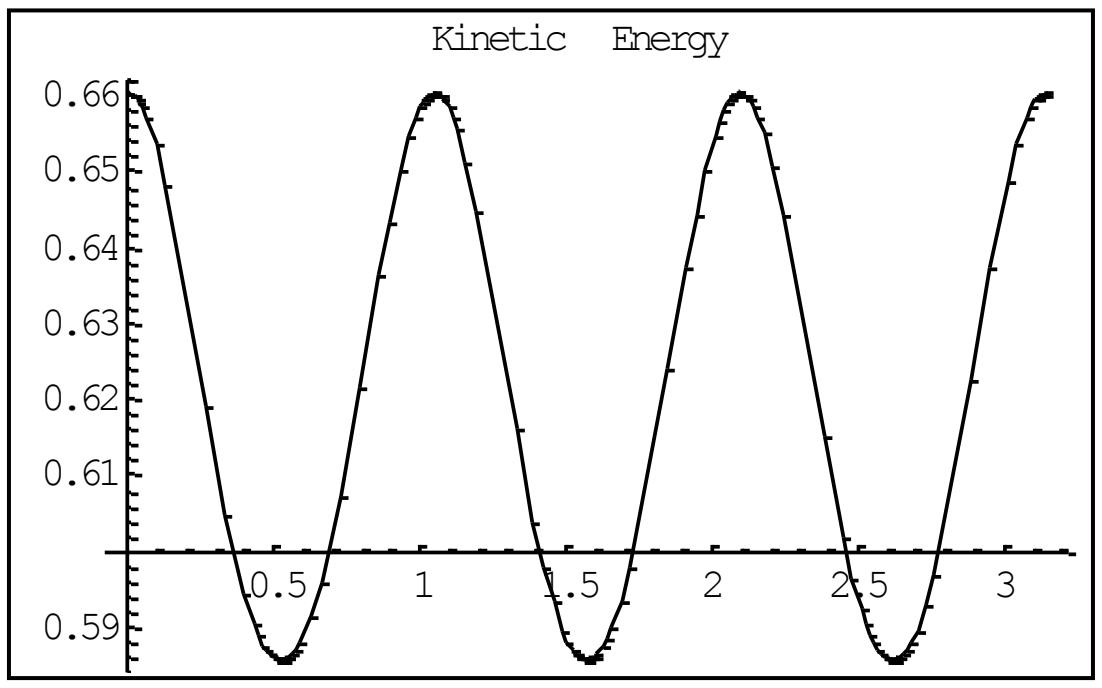

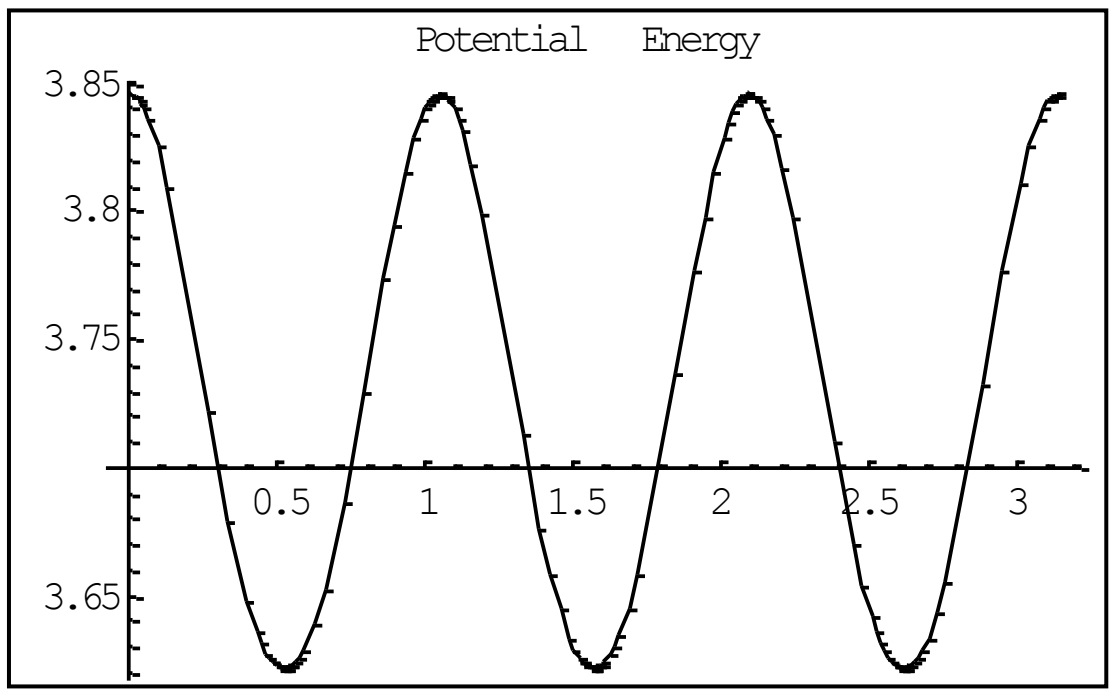

$$E_{k2r} = \frac{1}{2}V^2 ; E_{p2r} = \frac{G(m_1 + m_2)}{|r_2 - r_1|} . \quad (10)$$

**Fig.23 – Kinetic Energy of particle 1**

**Fig.24 – Potential Energy of the system**

The total energy is connected is then $E_2 = \frac{m_1 m_2}{m_1 + m_2}(E_{k2r} + E_{p2r})$ (11)

For two equal particles, $E_2 = M\, E_{T2r} /4$ with M=2m, the total mass. A similar relation does not exist in the three-body problem. The closest to it, is the relation between the COM of two of the particles and the third one.

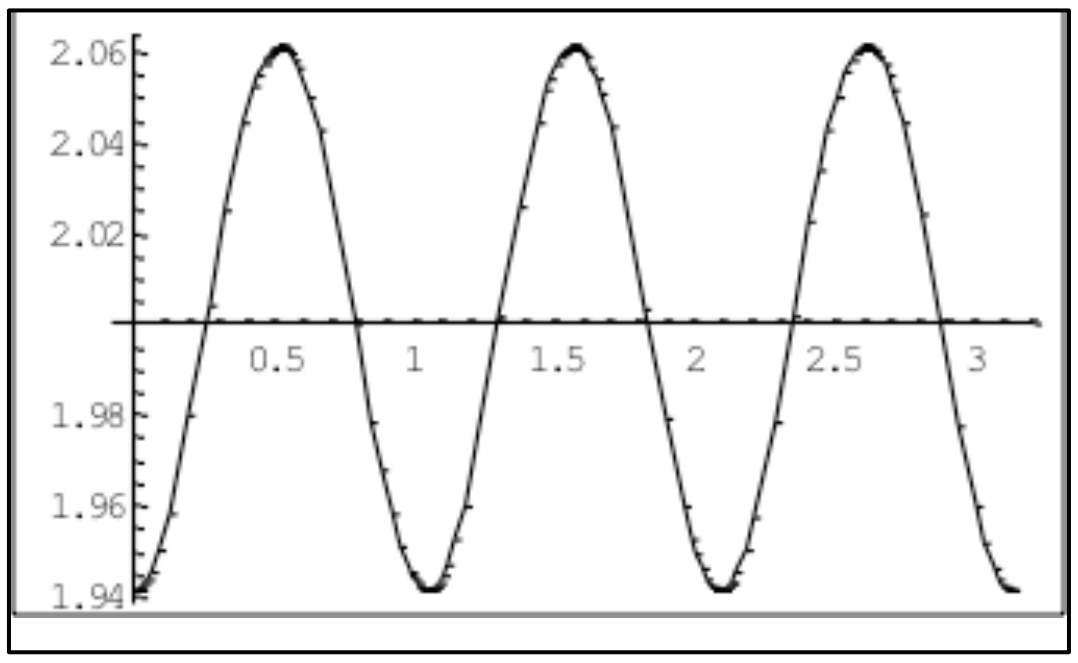

.

**Fig.25 – Ratio of Potential/Kinetic energy**

The ratio of the potential over the kinetic energy is shown in Fig.25.

P/K min =1.9413 1226 = 2 - .0586 8774

P/K max = 2.0612 4706 = 2 + .0612 4706

This ratio fluctuates around 2 which is the constant value for circular orbits in the two-body problem. For elliptic orbits, the fluctuations vary between 1/(1+e) and 1/(1-e).
In the 3BP, the well known collinear and triangular solutions (Lagrange points) have also the constant ratio of 2.

It is probably not possible to have trajectories with a constant value of this ratio when the velocities are not constant. When the fluctuations (Fig.25) are translated to an eccentricity, with the expressions from the two-body problem, it corresponds to e=.03. In the two-body problem, the eccentricity is usually interpreted as proportional (not linear) to the excess energy over the minimum energy required by the angular momentum. Here the angular momentum is zero. For the two-body problem, the eccentricity can also be interpreted as proportional to the exchange between kinetic and potential energy. The extrema occur at the collinear and isosceles configurations. In the two-body problem, the extrema occur at perigee and apogee.

*7.Time dependency*

Fig.26 shows the deviation of the perimeter of the triangle made up by the three points from a $\cos 6\tau$ function with adjusted parameters. Fig.27 shows similarly, the deviation from one side of the triangle from a $\cos 2\tau$ function.

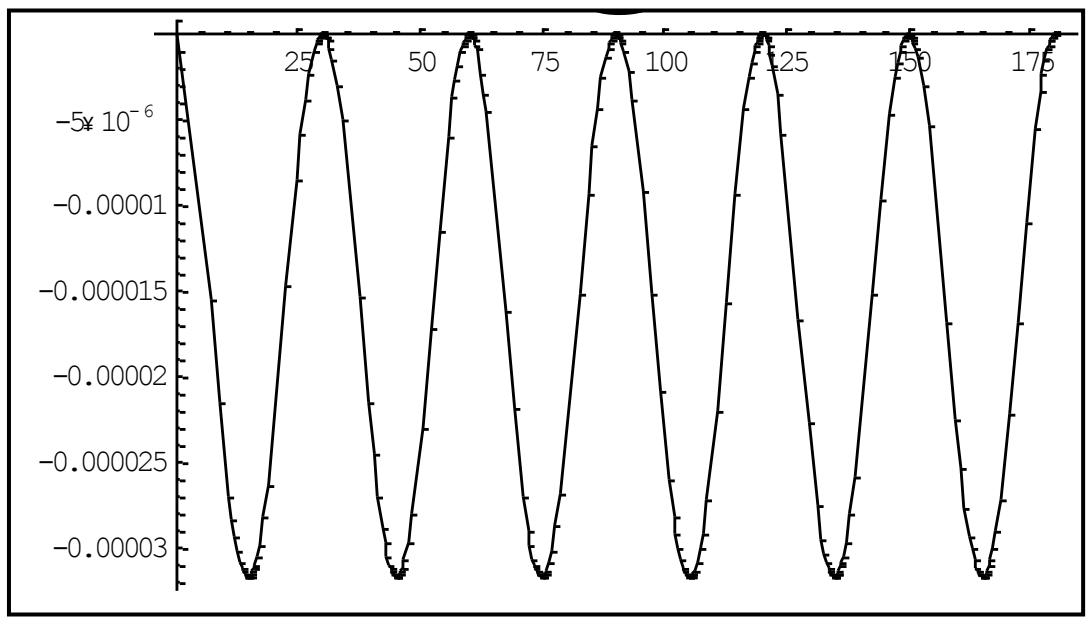

**Fig.21 – Deviation perimeter from $\cos 6\tau$**

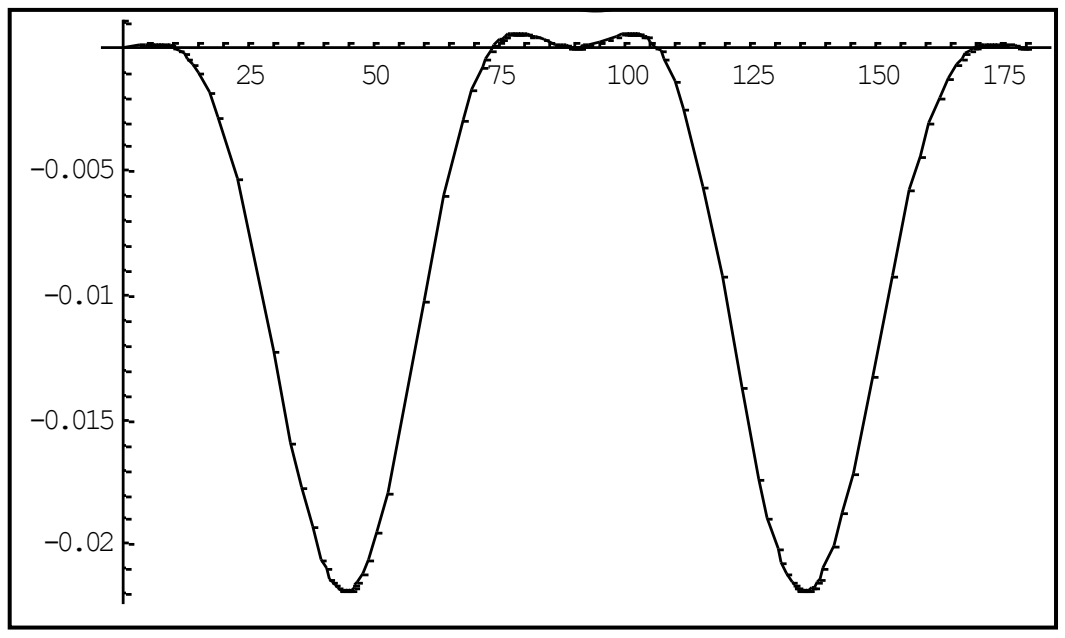

**Fig.22 – Deviation side 12 of the triangle**

## Summary

The figure-eight trajectory has the following intrinsic properties:

*1) Triangular configuration*

*The tangent at the point with x-coordinate = semi-major axis/2, meets the x- axis at the opposite extremal point of the curve..*

*The sides of the so-defined isoceles triangle, make an angle of 12.0202 deg with the x-axis*

*2) Collinear configuration*

*When one particle is at the origin, the three particles are on a line. The tangents at the two other points are parallel and have same direction as a tangents at the origin.*

*The direction of this line makes an angle of ±14.0688 deg with the x-axis.*

*The tangents at the origin make an angle of ±42.8434 deg with the x-axis.*

## Concluding Remarks

The discovery of the fig-8 solution of the three-body opened a new area of research in periodic solutions of the 3BP and nBP.

The fact that the three-body problem is not integrable does not exclude the possibility to describe some of these solutions by an appropriate family of curves. These curves will require a sufficient number of parameters that follow from the initial conditions. These parameters are not necessarily first integrals of the equations of motion. An illustrative example is the hypergeometric function that contains three free parameters $\alpha,\beta,\gamma$ . These parameters are not first integrals of the hypergeometric differential equation. The present situation may be compared to trying to find the solution to the two-body problem without any knowledge of conic sections except circles.

In this paper we collected some properties of the fig-8 trajectory and showed some possible lines of investigation.

References.

## Appendix A - Classical Figure-8 curves and the Fig8

The simplest way to get a fig-8 curve is to start from the parametric equation of an ellipse

$$\{x = aCosE \qquad y = bSinE\}$$

and double the frequency in the y-equation : $y = b\ Sin2E$

When b= a/2 , this curve is called *lemniscate of Gerono* and is also a *Lissjaous figure.* The triangular configuration with a base on x = a/2 (E=60) , has its top at x = - a / 4 and not at x = - a as in the fig-8 . This value is independent of the parameter b.

With the parameter b we can adjust the tangent at the origin. When a=2b, this angle is 45 deg . With b such that tan 42.84 = 2b/a or b = .46365 a , we have the angle of the fig-8 curve. Now, the collinear configuration occurs always for tan E = 1/ √2 => E = 35.264 and is also independent of b. This straight line makes an angle of 28.16 deg. which is practically twice the angle of the collinear line of the fig-8.

Table I gives various representations of well-known figure-8 curves. All three are of the fourth degree, have symmetry w.r.t. both axes and have a double point at the origin. A figure-8 curve that is not a fourth degree curve but defined by a complete elliptic integral of the first kind is :

$$\frac{K(k)}{\sqrt{(x+1)^2 + y^2}} = \frac{\pi}{2} \qquad k^2 = \frac{4x}{(x+1)^2 + y^2} \qquad \text{(A1)}$$

This curve shows up in the potential of a ring.

| | | Parametric Equation | Polar Equation | Cartesian Coordinates | Tangents O |
|---|---|---|---|---|---|
| L B | Lemniscate Bernoulli | $x = \frac{a\cos t}{1+\sin^2 t}$ <br> $y = \frac{a\cos t \sin t}{1+\sin^2 t}$ | $r^2 = a^2 \cos 2\theta$ | $(x^2 + y^2)^2 =$ <br> $a^2(x^2 - y^2)$ | $\pm 45 \deg$ |
| L G | Lemniscate of Gerono | $x = a\cos t_g$ <br> $y = a\cos t_g \ \sin t_g$ | $r^2 = a^2 \frac{\cos 2\theta}{\cos^4 \theta}$ | $x^4 = a^2(x^2 - y^2)$ | ±45 deg |
| H P | Hippopede of Proclus (k > 1) | $x = a\sqrt{1 - k^2 \sin^2 t}\ \cos t$ <br> $y = a\sqrt{1 - k^2 \sin^2 t}\ \sin t$ | $r^2 = a^2(1 - k^2 \sin^2 t)$ | $(x^2 + y^2)^2 =$ <br> $a^2(x^2 - (k^2 - 1)y^2$ | ±atan $\frac{1}{\sqrt{k^2 - 1}}$ |

**Table I. – Some figure-eight curves.**

The scaling parameter a is the maximum value on the x-axis, so a = .746156. Working with a = 1 is convenient for the analysis of the properties of those curves. The reference length R is then 1.3402 a.

The lemniscate of Bernoulli (LB) and the lemniscate of Gerono (LG) have only “a” as free parameter and their tangents at the origin are fixed to ±45 deg.

### ***3.1 Simple example of an adjustment with hippopede Proclus (HP)***

As the hippopede has a second parameter k, we can easily make the tangents at the origin fit with the fig8 of the 3BP. With k = 1 / (sin 42.8434) = 1.4706, the tangents at the origin have the correct value (Fig.A1).When $k^2$= 2, HP reduces LB which is evident from its polar and Cartesian equation. In the parametric representation, the meaning of the parameter t is different.

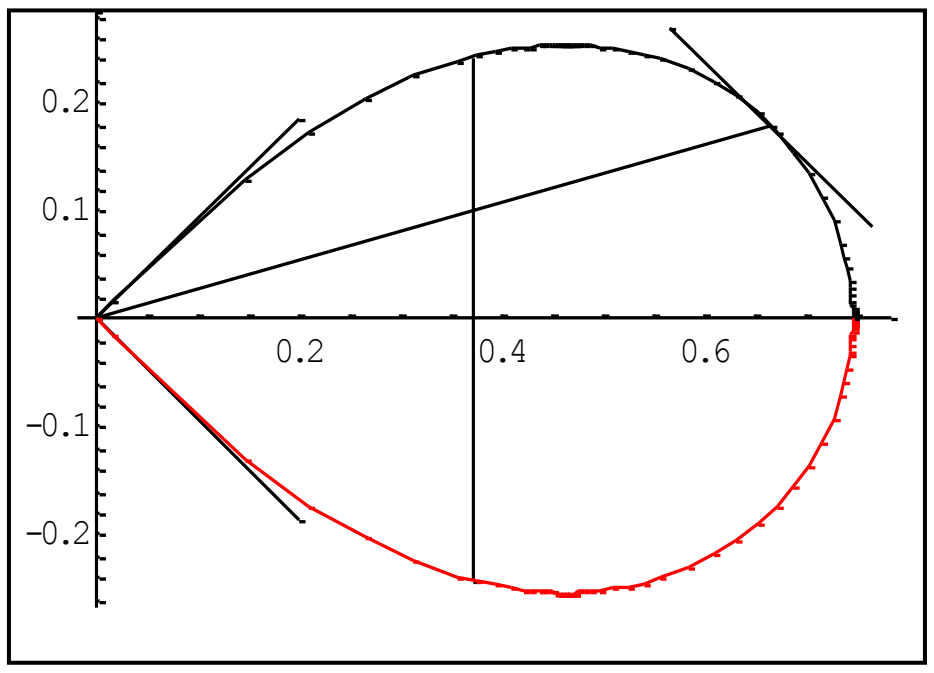

**Fig.A1 –Hippopede with adjusted k**

Table II compares the other relevant values of the adjusted HP with those of the fig-8

| | Collinear | | Triangular | | |
|---|---|---|---|---|---|
| | Tangent at O | Point with same tangent | Base point | tangent | height |
| 3BP | 42.8434 | 14.0688<br>r =.690236 | x= a/2<br>y=.238393 | 12.02 | 1.5 a |
| HP | 42.8434 | 14.9532<br>r=.690351 | x= a/2<br>y=.24294 | 13.2732 | 1.38025 a |

**Table II. Comparison adjusted hippopede k=1.4706 and three-body values**

The Collinear configuration makes an angle of 14.9532 deg. with the x-axis as compared to the 14.0688 deg. For the triangular configuration, the height of the isosceles triangle is 1.38 a instead of the required 1.5 a.

A new free parameter can be added by changing a to a' on the y-coordinate only.
The parameter a', changes the equation of the hippopede to:

$$(\rho^2 x^2 + y^2)^2 = a'^2 (\rho^2 x^2 - (k^2 - 1) y^2) \quad \text{with } \rho = a'/a \qquad \text{(A2)}$$

Now, the direction of the collinear configuration can be adjusted to 14.0688. Table III shows the improvement in the approximation to the fig-8 trajectory. The value of the parameter t in the hippopede is close to the values of the polar angle θ in the collinear and triangular configuration: $t_c$=14.994 and $t_T$=34.134.

| | Collinear | | Triangular | | |
|---|---|---|---|---|---|
| | Tangent at O | Point with same tangent | Base point | tangent | meets x-axis at |
| HP | 42.8434 | 14.0688<br>r=.691006 | x= a/2<br>y=.2365 | 12.0473 | 1.4855 a |

**Table III – Hippopede adjusted with k=1.4204, a'=.697886**

## *3.2 Lemniscate of Bernoulli (LB)*

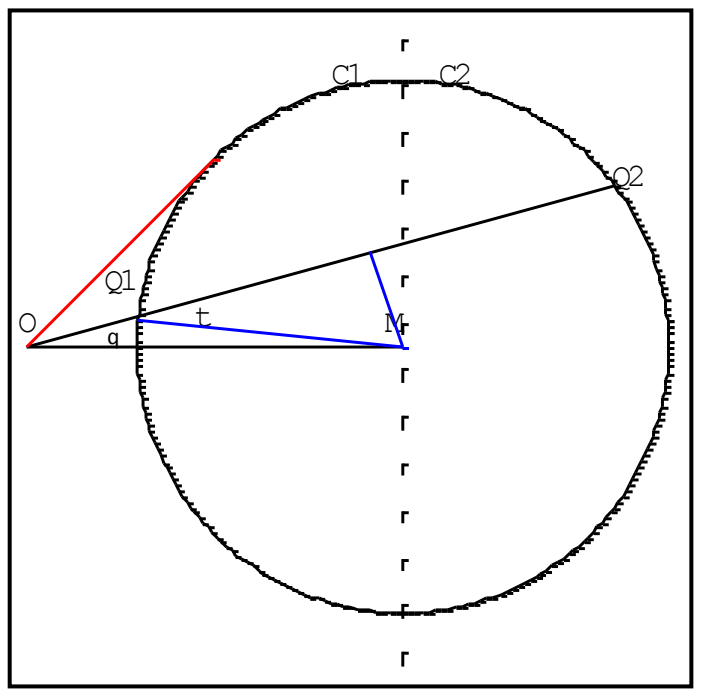

**Fig.A2 – Cissoid construction**

The LB can be constructed as a "cissoid" : from O w.r.t. two given curves each of which gives an intersection point $Q_1$ and $Q_2$. In this case the two curves are the two half-circles ($C_1$,$C_2$ ) that make up the given circle. A line from O cuts the circle in 2 points $Q_1$ and $Q_2$. Laying off the segment $Q_1Q_2$ from O in its direction defines the cissoid.
The maximum value of this segment is the diameter of the circle $a = 2r_a$ and is the semi-major axis of the curve.
The maximum value $\theta_m$ of the angle $\theta$ that gives two points $Q_1$ and $Q_2$ is $\underline{\sin\theta_m = 1/k}$ degrees. $\theta_m$ is also the angle of the tangents at the origin.

When $k^2$=2, the curve is LB and $\theta_m$ = 45 deg.
From this definition one obtains easily the polar representation via a parametric representation in t

$$r = (Q_1Q_2) = a\cos t' \quad \sin\theta = \frac{\sin t'}{k} \quad \Rightarrow r = a\sqrt{1-k^2\sin^2\theta} \qquad \text{(A3)}$$

The parameter t' varies from 0 to π/2 when one fourth of the curve is described. The Cartesian equations are :

$$\begin{aligned} x = r\cos\theta \quad &= \frac{a}{k}\sqrt{k^2-\sin^2 t'}\cos t' \\ y = r\sin\theta \quad &= \frac{a}{k}\cos t'\sin t' \\ \mathit{Cartesian} \Rightarrow \quad & (x^2+y^2)^2 = a^2(x^2-(k^2-1)y^2) \end{aligned} \qquad \text{(A4)}$$

Table I gives the standard parametric form [Lawrence][2] where this curve is defined as the inverse of an equilateral hyperbola. . In the parameterization t' , the LB is :

$$r = a\cos t \quad \sin\theta = \frac{\sin t}{\sqrt{2}} \qquad \text{(A4b)}$$

The parameters t and t' are related by :

$$\tan t = \frac{\tan t'}{\sqrt{k^2+(k^2-2)\tan^2 t'}} \qquad \text{(A5)}$$

which reduces to tan t = 1/√2 tan t' for LB.

The tangents at the origin follow also from the limit of the polar angle as the lemniscate degenerates to a straight line at the origin:

$$\lim_{t'\to\pi/2} \tan\theta = \frac{\sin t'}{\sqrt{1+\cos^2 t'}}\ 1 \ \Rightarrow\ \theta = \pi/4 \qquad \text{(A6)}$$

When k=1/sin 42.8434=1.4706 we have a figure-8 curve with tangents at the origin as the fig-8.

The direction of the tangent at an arbitrary point of LB in terms of the parameter t and the polar angle θ, is :

$$y' = -\frac{1}{\sin t}\frac{1-3\sin^2 t}{3-\sin^2 t} = -\frac{1}{\tan\theta}\frac{1-3\tan^2\theta}{3-\tan^2\theta} = -\frac{1}{\tan 3\theta} \qquad \text{(A7)}$$

This formula shows that the normal to the tangent makes always an angle 3θ with the x-axis.

Hence, any tangent makes an angle 3θ+π/2 with the x-axis. The polar angle of the collinear configuration is π/2+3θ = is π/4+π/2 => 15 deg.
On the fig-8 $\theta_c$ =14.0688 deg. For a cissoid with k

$$\sin 2\theta_c = \frac{\sin 2\theta_M}{2} = \frac{\sqrt{k^2-1}}{k^2} \quad as \quad \sin\theta_M = \frac{1}{k}$$ is the angle at the origin.

k=1.22525, is required the collinear polar angle of the fig-8. This value differs from 1.4706 that is required for the fig-8.

**In passing we notice also that for any angle tan3x = - tan(x-60) tan x tan (x+60).**

1) *Triangular property*

The point of LB with *x-coordinate 1/2, (a=1) has y-coordinate* y = tan $\theta_T$ /2
The Cartesian Equation (A4) shows that , $\theta_T$ is the solution of :

$$\cos^4\theta_T - \tfrac{1}{4}\cos^2\theta_T - 1/8 = 0 \qquad \Rightarrow \cos^2\theta_T = (1+\sqrt{3})/4 \Rightarrow \theta_T = 34.26465 \quad \text{(A9)}$$

Using this result in Eq.A7, combined with the fact that the tangent y' makes an angle 3θ- π/2 with the axis , the segment H= [x=1/2, point where the tangent cuts the axis ] is easily computed. The result is H=3/2 . This is the point x=-1 and *the LB has the triangular property.*
(on the fig-8, $\theta_T$ =32.569)

2) *The pedal equation of LB is $p = r^3$*

Where p is the perpendicular from the origin on a tangent. The segment p is parallel to the normal at the point of the tangent line and has a meaning in the angular momentum h : $h = r^2\dot{\theta} = r\,v\cos\ \tau = p\,v\ and\ \frac{p}{r} = \cos\tau\ = \cos 2\theta$

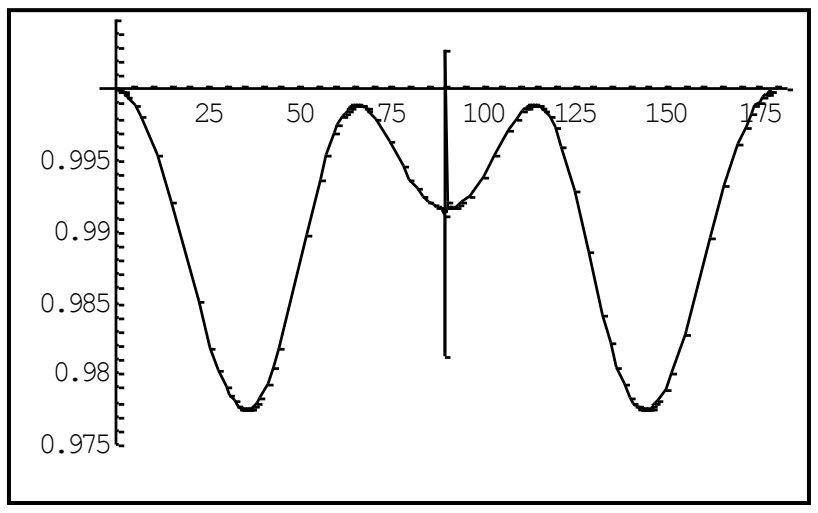

**Fig.A3 – Pedal equation LB applied to fig-8**

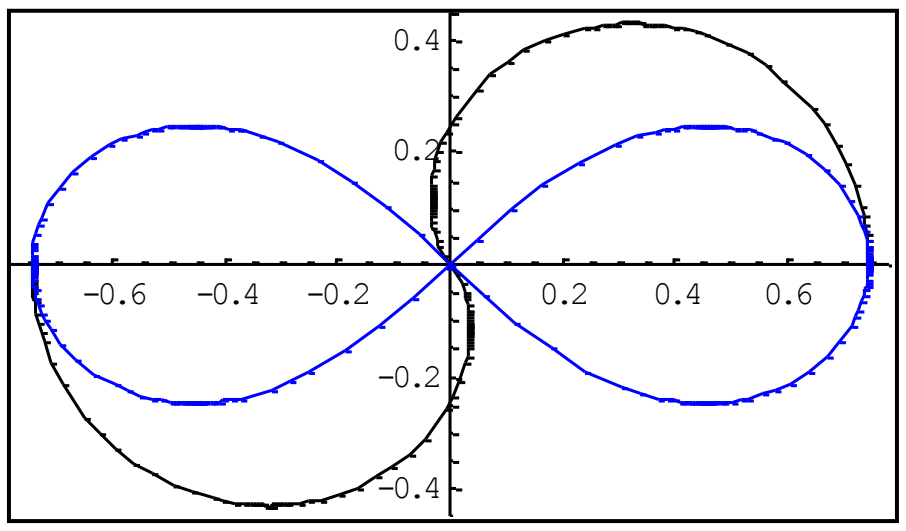

**Fig.A4 – Pedal Curve of LB**

Fig.A2 shows the ratio p/$r^3$ of the fig-8. This ratio fluctuates between 1 and .9775 at t = 35.72. There is also a relative maximum of .99893 at t = 65.8702.
Fig.A3 shows the pedal curve.

### *3) Arclength and Elliptic Integrals*

The arclength of LB is given by an elliptic integral of the first kind with modular angle 45 deg or parameter ½.:

$$ds = \sqrt{1 + r^2\left(\frac{d\theta}{dr}\right)^2}\,dr = \frac{dr}{\sqrt{1-r^4}} \quad as \quad \frac{d\theta}{dr} = -\frac{r}{\sqrt{1-r^4}} \qquad \text{(A10)}$$

So $L_{lb}$ = 2 a √2 K[1/√2] = 2.622057 K = complete elliptic integral of the first kind.

For the fig8 $L_{f8}$ = 2a 2.55968 17440

From the integration, we can easily construct the curve s(t). and compare it toEq.( A10).

$$ds = \frac{v\, r}{x v_x + y v_y}\, dr \qquad \text{(A11)}$$

When the radical $\sqrt{1-r^4}$ is replaced by $\sqrt{(1-r^2)(1+mr^2)}$ either for dθ/dr or for ds/dr, we define curves with a parameter m whose arclength remain elliptic integrals. For m=1 we recover LB.

As a first application we used such radical for dθ/dr . The corresponding ds/dr is :

$$\frac{ds}{dr} = \sqrt{\frac{1-(1-m)r^2(1-r^2)}{(1-r^2)(1+mr^2)}} \qquad \text{(A12)}$$

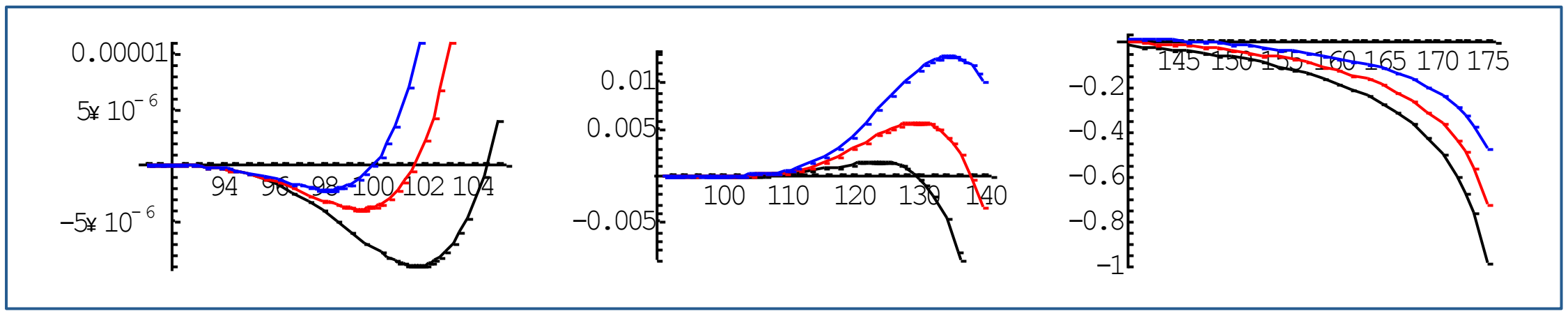

**Fig.A5 – Comparison ds/dr_calculated – ds/dr_assumed for m=1.1,1.2,1.3 (eq.S3)**

Fig.A5 shows zooms of the difference between ds/dr from the trajectory and ds/dr as assumed by eq.S3. The values of m are choosen such as to minimize the deviations. To start with r = 0, the results of the integration from π/2 to π are used. On the global curves, the differences are hardly noticeable. Fig.A5 shows that the deviation is first negative (τ = 90, ~ 102) , then positive (τ~103, ~140) with a larger excursion and then going negative with an increasing value when τ approaches 180. As we approach 180 degrees, r reaches an extremum (dr = 0) , dθ/dr does not make much sense and the investigation much be reformulated.

### *4) Comparison fig8 and LB*

In Table III parameter t refers to the presentation in Table I and parameter t' to the presentation in eq. (A3). From the integration we have also the “time” of a particular event normalized as an angle. The corresponding values of the parameters t and t 'are different. However the parameter $t_g$, appearing in the lemniscate of Gerono (LG), has respectively the values 30,60,90 for the successive collinear and triangular configurations although the points themselves are quite different. This suggests that a parameterization of {x,y} by $t_g$ , is possibly related to the “mean anomaly” , one of the other angles (t',p) can be related to an “eccentric anomaly” while θ plays the role of “true anomaly”.

| | θ (deg.)LB | θ (deg.)3B | ratio r | param t | param p | time τ | $t_g$ (LG) |
|---|---|---|---|---|---|---|---|
| Collinear<br>+tangent | 15<br>45 | 14.0688<br>42.8434 | 1.066189<br>1.050337 | 21.4707 | 15.542 | 30 | 30 |
| Min $p/r^3$ | - | 20.7048 | | 30 | 22.2076 | 35.7199 | |
| | | | | | | | |
| extrem. y value (a=1) | 30<br>.353553 | 28.5123 | 1.052177 | 45 | 35.2644 | 53.6628 | 45 |
| Triangular<br>+tangent | 34.2647<br>12.7939 | 32.569<br>12.02 | 1.052065<br>1.064384 | 52.7711 | 42.941 | 60 | 60 |
| Max $p/r^3$ | - | 37.7612 | | 60 | 50.768 | 65.87 | |
| | | | | | | | |
| Length a=1 | 2.622057 | 2.559682 | 1.024368 | | | | |
| tan-O/tang-C | 45/15<br>=3 | 42.84/14.068<br>= 3.045277 | 1.015092 | | | | |

**Table III . – Lemniscate of Bernoulli and 3body curve**

The value of $t_g$ for the extremum of y is 45 deg. On the trajectory, this extremum occurs at τ = 53.6628 which shows that a simple identification of $t_g$ with τ does not work.

Projecting the LB in a plane rotated by an angle α about the x-axis, the x-coordinates remain

| | r1 | r2 | r3 | $\Sigma r_i^2$ | $E_p$ | $E_k$ | \|P/K\| |
|---|---|---|---|---|---|---|---|
| T | 1 | .593303 | .593303 | 1.704017 | -1.28195 | .66305 | 1.94131 = 2 - .05869 |
| C | 0 | .925055 | .925055 | 1.711453 | -1.20732 | .585722 | 2.06125 = 2 + .06125 |
| | | | | | | | |

unchanged and the y-coordinates are reduced by k’= cos α . Choosing an appropriate a, the differences between this curve and the Fif-8 can be reduced to some extend.